\documentclass{svjour3}                     
\smartqed  
\usepackage{graphicx}
\usepackage{amssymb,amsmath,amsfonts,latexsym,euscript}
\usepackage{xcolor}

\usepackage{multirow}
\usepackage{algorithm}
\usepackage{algpseudocode}
\usepackage{float}
\newfloat{algorithm}{tbp}{lap}

\def\d{\,{\rm d}}
\def\E{\,{\rm e}}

\def\i{{\rm i}}

\begin{document}

\title{{Optimal quadrature formulas for non-periodic functions in Sobolev space and its application} {to CT image reconstruction}}
\thanks{A.R. Hayotov is supported by ISEF program at Korea Foundation for Advanced Studies.}

\titlerunning{Optimal quadrature formulas for CT image reconstruction}        

\author{Abdullo R. Hayotov$^*$, Soomin Jeon, Chang-Ock Lee, Kholmat M. Shadimetov}


\institute{Abdullo R. Hayotov \at
$^*$Corresponding author:\\
Department of Mathematical Sciences, KAIST, 291 Daehak-ro, Yuseong-gu, Daejeon 34141, Republic of Korea,\\
V.I.Romanovskiy Institute of Mathematics, Uzbekistan Academy of Sciences,
81 M.Ulugbek str., Tashkent 100170, Uzbekistan,\\
\email{hayotov@mail.ru}.
\and
Soomin Jeon \at
Department of Mathematical Sciences, KAIST, 291 Daehak-ro, Yuseong-gu, Daejeon 34141, Republic of Korea \\
\email{soominjeon@kaist.ac.kr}.
\and
Chang-Ock Lee \at
Department of Mathematical Sciences, KAIST, 291 Daehak-ro, Yuseong-gu, Daejeon 34141, Republic of Korea \\
\email{colee@kaist.edu}.
\and
Kholmat M.Shadimetov \at
V.I.Romanovskiy Institute of Mathematics, Uzbekistan Academy of Sciences,
81 M.Ulugbek str., Tashkent 100170, Uzbekistan,\\
Tashkent Railway Engineering Institute, 1 Odilxojaev str., Tashkent 100167, Uzbekistan\\
\email{kholmatshadimetov@mail.ru}.
}

\date{Received: date / Accepted: date}

\maketitle
\begin{abstract}
In the present paper, optimal quadrature formulas in the sense of Sard are constructed for numerical integration of the integral $\int_a^b\E^{2\pi\i\omega x}\varphi(x)\d x$ with $\omega\in \mathbb{R}$ in the Sobolev space $L_2^{(m)}[a,b]$ of complex-valued functions which are square integrable with $m$-th order derivative. Here, using the discrete analogue of the differential operator $\frac{\d^{2m}}{\d x^{2m}}$, the explicit formulas for optimal coefficients are obtained. The order of convergence of the obtained optimal quadrature formula is $O(h^m)$.
{As an application, we implement the filtered back-projection (FBP) algorithm, which is a well-known image reconstruction algorithm for computed tomography (CT).
By approximating Fourier transforms and its inversion using the proposed optimal quadrature formula} {of the second and third orders} {, we observe that the accuracy of the reconstruction algorithm is improved.
In numerical experiments, we compare the quality of the reconstructed image obtained by using the proposed optimal quadrature formulas with the conventional FBP, in which} {fast Fourier transform is used} {for the calculation of Fourier transform and its inversion.
In the noise test, the proposed algorithm provides more reliable results against the noise than the conventional FBP.}


\end{abstract}

\maketitle

\section{Introduction}

It is known that when complete continuous X-ray data are available Computed Tomography (CT) image can be reconstructed exactly using the filtered back-projection formula (see, for instance, \cite{Buzug08,Feeman15,KakSlaney88,Nat2001}). This formula gives interactions between the Radon transform, the Fourier transform and the back-projection transform. A description of the filtered back-projection formula along \cite[Chapter 3]{KakSlaney88} is provided below.

In the Cartesian system with $x,y$-axes consider a unit vector $(\cos\theta,\sin\theta)$. Then the line perpendicular to this vector with the distance $t$ to the origin can be expressed as ${{\ell }_{t,\theta }}$: $x\cos \theta +y\sin \theta =t$. Assume the object is represented by a two variable function $\mu(x,y)$, which denotes the attenuation coefficient in X-ray CT applications. Then, the $\theta $-view projection along the line ${{\ell }_{t,\theta }}$ can be expressed as
$$
P(t,\theta )=\int\limits_{-\infty }^{\infty }{\int\limits_{-\infty }^{\infty }{\mu(x,y)\delta (x\cos \theta +y\sin \theta -t)\d x\d y}},
$$
where $\delta$ denotes the Dirac delta-function. The function $P(t,\theta )$ is known as the Radon transform of $\mu(x,y)$. A projection is formed by combining a set of line integrals. The simplest projection is a collection of parallel ray integrals as is given by $P(t,\theta )$ for a constant $\theta $. This is known as a parallel beam projection. It should be noted that there are fan-beam in 2D and cone-beam in 3D projections \cite{Buzug08,Feeman15,KakSlaney88}.
	
The problem of CT is to reconstruct the function $\mu(x,y)$ from its  projections $P(t,\theta )$. There are analytic and iterative methods for CT reconstruction. One of the widely used analytic methods of CT reconstruction is the filtered back-projection method. It can be modeled by
\begin{equation}\label{eq1.1}
\mu(x,y)=\int\limits_{0}^{\pi }{\int\limits_{-\infty }^{\infty }{S(\omega ,\theta )\left| \omega  \right|{{\E}^{2\pi \i\omega
(x\cos \theta +y\sin \theta )}}\d\omega \d\theta ,}}
\end{equation}
where
\begin{equation}\label{eq1.2}
S(\omega ,\theta )=\int\limits_{-\infty }^{\infty }{P(t,\theta ){{\E}^{-2\pi \i\omega t}}\d t}		
\end{equation}
is the 1D Fourier transform of $P(t,\theta )$. The inner integral of (\ref{eq1.1}),
\begin{equation}\label{eq1.3}
Q(t,\theta )=\int\limits_{-\infty }^{\infty }{S(\omega ,\theta )|\omega|{{\E}^{2\pi \i\omega t}}\d\omega},
\end{equation}		
is a 1D inverse Fourier transform of the product $S(\omega ,\theta )\left| \omega  \right|$,
which represents a projection filtered by a 1D filter whose frequency representation is $|\omega|$. The outer integral performs back-projection. Therefore, the filtered back-projection consists of two steps: filtration and then back-projection.

Thus,  the Fourier transforms play the main role in (\ref{eq1.1})-(\ref{eq1.3}). But in practice, due to the fact that we have discrete values of the Radon transform, we have to approximately calculate the Fourier transforms in the filtered back-projection.
For this purpose, it is necessary to consider the problem of approximate calculation of the integral
\begin{equation}\label{eq1.4}
I(\varphi)=\int\limits_a^b\E^{2\pi \i\omega x}\varphi(x)\d x
\end{equation}
with $\omega\in \mathbb{R}$. This type of integrals are called \emph{highly oscillating integrals}.
In most cases it is impossible to get the exact values of such integrals.
Thus, they can be approximately calculated using the formulas of numerical integration.
However, standard methods of numerical integration cannot be successfully applied for that.
Therefore special effective methods should be developed for approximation of highly oscillating integrals.
One of the first numerical integration formula for the integral (\ref{eq1.4}) was obtained by Filon \cite{Filon28} in 1928 using a quadratic spline.
Since then, for integrals of different types of highly oscillating functions many special effective methods have been developed, such as
Filon-type method, Clenshaw-Curtis-Filon type method, Levin type methods, modified Clenshaw-Curtis method, generalized quad\-rature rule, and Gauss-Laguerre quadrature (see, for example, \cite{AvdMal89,IBab,BabVitPrag69,BakhVas68,IserNor05,Mil98,NovUllWoz15,XuMilXiang,ZhangNovak19}, for more review see, for instance, \cite{DeanoHuyIser18,MilStan14,Olver08} and references therein).
Recently, in \cite{BolHayMilShad17,BolHayShad16,BolHayShad17}, based on Sobolev's method, the problem of construction of optimal quadrature formulas
in the sense of Sard for numerical calculation of integrals (\ref{eq1.4}) with integer $\omega$ was studied in Hilbert spaces $L_2^{(m)}$ and $W_2^{(m,m-1)}$.
Here, we consider the Sobolev space $L_2^{(m)}[a,b]$ of non-periodic, complex-valued functions defined on the interval
$[a,b]$, which possess an absolute continuous $(m-1)${-th} derivative on $[a,b]$, and whose $m${-}th order derivative is square integrable \cite{Sobolev74,SobVas}.
The space $L_2^{(m)}[a,b]$ is a Hilbert space with the inner product
\begin{equation}\label{eq1.5}
\langle\varphi,\psi\rangle=\int\limits_a^b\varphi^{(m)}(x)\ \bar\psi^{(m)}(x)\d x,
\end{equation}
where $\varphi^{(m)}$ is the $m$-th order derivative of the function $\varphi$ with respect to $x$, $\bar\psi$ is the complex conjugate function to the function $\psi$ and the norm of the function $\varphi$ is correspondingly defined by the
formula
$$
\|\varphi\|_{L_2^{(m)}[a,b]}=\langle\varphi,\varphi\rangle^{1/2}.
$$

The aim of the present work is to construct optimal quadrature formulas in the sense of Sard in the Sobolev space $L_2^{(m)}$ for numerical integration of the integral (\ref{eq1.4}) with real $\omega$ using the results of the work \cite{BolHayShad17}.
Then to approximate reconstruction of CT image, we apply the obtained optimal quadrature formulas
to the Fourier transform and its inversion in (\ref{eq1.1}).
Recently we got the results  when $m=1$ \cite{HayJeonLee19}. There is a work by Reider and Faridani \cite{Faridan03} for an optimal reconstruction algorithm, where the optimal $L^{2}-$convergence rates of filtered back-projection (FBP) algorithm was provided with assumption that the Radon transform was analytically computable.

The rest of the paper is organized as follows. Section 2 is devoted to construction of optimal quadrature formulas in the sense
of Sard in the space $L_2^{(m)}$ for numerical calculation of Fourier integrals. There are obtained analytic formulas for optimal coefficients
using the discrete analogue of the differential operator $\frac{\d^{2m}}{\d x^{2m}}$.
In section 3, the obtained optimal quadrature formulas for the cases $m=2$ and $m=3$ are applied for CT image reconstruction by approximating
Fourier transforms in the filtered back-projection formula.

\section{Optimal quadrature formulas for Fourier integrals in the space $L_2^{(m)}[a,b]$}
\setcounter{equation}{0}

In $L_2^{(m)}[a,b]$ space for approximation of the integral (\ref{eq1.4}), we consider the following quadrature formula
\begin{equation}\label{eq2.1}
\int\limits_a^b \E^{2\pi \i\omega x}\varphi (x)\d x \cong  \sum\limits_{\beta  = 0}^N
{C_\beta} \varphi (x_\beta)
\end{equation}
with the error
\begin{equation}\label{eq2.2}
(\ell,\varphi) = \int\limits_a^b {\E^{2\pi\i \omega  x}\varphi (x)\d x -
\sum\limits_{\beta  = 0}^N {C_\beta  \varphi (x_\beta)}},
\end{equation}
where
$
(\ell,\varphi)=\int\limits_{-\infty}^\infty  {\ell (x)\varphi
(x)\d x}
$ is the value of the error functional $\ell$ at the function $\varphi$. Here the error functional $\ell$ has the form
\begin{equation}\label{eq2.3}
\ell (x) =\E^{2\pi \i\omega x}\varepsilon _{[a,b]} (x) - \sum\limits_{\beta  = 0}^N
{C_\beta} \delta (x - x_\beta),
\end{equation}
$C_\beta$ are coefficients, $x_{\beta}=h\beta+a$ ($\in [a,b]$) are nodes of the formula (\ref{eq2.1}), $h=\frac{b-a}{N}$,
$N\in \mathbb{N}$, $\i^2=-1$, $\omega\in \mathbb{R}$,  $\varepsilon
_{[a,b]}(x)$ is the characteristic function of the interval $[a,b]$, and $\delta$ is the Dirac delta-function.
We mention that the coefficients $C_{\beta}$ depend on $\omega$, $h$ and $m$.

The error of the formula (\ref{eq2.1}) defines a linear functional in $L_2^{(m)*}[a,b]$, where $L_2^{(m)*}[a,b]$ is the
conjugate space to the space $L_2^{(m)}[a,b]$. Since the functional (\ref{eq2.3}) is defined on the space $L_2^{(m)}$, the
conditions
\begin{equation}\label{eq2.3-1}
(\ell,x^{\alpha})=0,\ \alpha=0,1,...,m-1
\end{equation}
should be fulfilled. The conditions (\ref{eq2.3-1}) mean the exactness of the quadrature formula (\ref{eq2.1}) for all algebraic polynomials
of degree less than or equal to $m-1$. Hence we get that for a function from the space $C^m[a,b]$ the order of convergence
of the optimal quadrature formula (\ref{eq2.1}) is $O(h^m)$.

It should be mentioned that from equalities (\ref{eq2.3-1}) one can get  the condition $N+1\geq m$ for existence
of the optimal quadrature formula of the form (\ref{eq2.1}) in the space $L_2^{(m)}$.

Sard's optimization problem of numerical integration formulas of the form (\ref{eq2.1}) in the space $L_2^{(m)}[a,b]$
is the problem of finding the minimum of the norm of the error functional $\ell$ by coefficients $C_{\beta}$, i.e.,
to find coefficients $C_{\beta}$ satisfying the equality
\begin{equation}\label{eq2.3-2}
\|\mathring{\ell}\|_{L_2^{(m)*}[a,b]}=\inf\limits_{C_\beta}\|\ell\|_{L_2^{(m)*}[a,b]}.
\end{equation}
The coefficients satisfying the last equality are called \emph{optimal coefficients} and they are denoted as $\mathring{C}_{\beta}$.
The  quadrature formula with coefficients $\mathring{C}_{\beta}$ is called \emph{the optimal quadrature formula in the sense of Sard}, and
$\mathring{\ell}$ is the error functional corresponding to the optimal quadrature formula.

The solution of Sard's problem gives the sharp upper bound for the error (\ref{eq2.2}) of functions $\varphi$ from the space $L_2^{(m)}[a,b]$ as follows
$$
|(\ell,\varphi)|\leq \|\varphi\|_{L_2^{(m)}[a,b]}
\|\mathring{\ell}\|_{L_2^{(m)*}[a,b]}.
$$

This problem, for the quadrature formulas of the form (\ref{eq2.1}) with $\omega=0$, was first studied by Sard
\cite{Sard} in the space $L_2^{(m)}$ for some $m$. Since then, it was investigated by many authors (see, for instance, \cite{CatCom,Com72a,Com72b,GhOs,FLan,MalOrl,ShadHay11}) using spline method, $\phi$-function method, and Sobolev's method. Finally, in the works \cite{Koh,Shad83,Shad99} this problem was solved for any $m\in \mathbb{N}$ with equally spaced nodes and the explicit expressions for the optimal coefficients have been obtained (see Theorem \ref{Thm4}).

It should be noted that the problem of construction of lattice optimal cubature formulas in the space $L_2^{(m)}$ of multi-variable functions was first stated and investigated by Sobolev \cite{Sobolev74,SobVas}. Further, in this section, based on the results of the work \cite{BolHayShad17}, we solve Sard's problem on construction of
optimal quadrature formulas of the form (\ref{eq2.1}) for $\omega\in \mathbb{R} $ with $\omega\neq 0$,
first for the interval $[0,1]$ and then using a linear transformation for the interval $[a,b]$.
For this we use the following auxiliary results.

\subsection{Preliminaries}

We need the concept of discrete argument functions and operations on them.
For this we give the definition for functions of discrete argument by following  \cite{Sobolev74,SobVas}.

Assume that the nodes $x_\beta$ are equally spaced, i.e.,
$x_\beta=h\beta,$ $h$ is a small positive parameter, and $\varphi$ and
$\psi$ are complex-valued functions defined on the real line $\mathbb{R}$ or on an interval of $\mathbb{R}$.
The function $\varphi (h\beta )$ given on some set
of integer values of $\beta$ is called {\it a function of discrete argument}.
{\it The inner product} of two discrete argument
functions $\varphi(h\beta )$ and $\psi (h\beta )$ is defined by
$$
\left[ {\varphi(h\beta),\psi(h\beta) } \right] =
\sum\limits_{\beta  =  - \infty }^\infty  {\varphi (h\beta ) \cdot
\bar \psi (h\beta )}
$$
{if the series on the right hand side of the last equality
converges absolutely.}
{\it The convolution} {of two
functions $\varphi(h\beta )$ and $\psi (h\beta )$ is the following inner
product}
$$
\varphi (h\beta )*\psi (h\beta ) = \left[ {\varphi (h\gamma ),\psi
(h\beta  - h\gamma )} \right] = \sum\limits_{\gamma  =  - \infty
}^\infty  {\varphi (h\gamma ) \cdot \bar\psi (h\beta  - h\gamma )}.
$$

We note that coefficients of optimal quadrature formulas and interpolation splines
in the spaces $L_2^{(m)}$ and $W_2^{(m,m-1)}$ depend on the roots of the Euler-Frobenius type polynomials
(see, for instance, \cite{BabHay19,BolHayShad17,CabHayShad14,Koh,Shad10,ShadHay11,ShadHay13,ShadHay14,ShadHayAkhm15,ShadHayNur13,ShadHayNur16,Sobolev06,SobVas}).
\emph{The Euler-Frobenius polynomials} $E_k (x)$, $k = 1,2,...,$ are
defined as follows (see, for instance, \cite{Frob1910,SobVas}):
\begin{equation}\label{eq2.5}
 E_k (x) = \frac{{(1 - x)^{k + 2} }}{x}\left( {x\frac{\d}{{\d x}}}
\right)^k \frac{x}{{(1 - x)^2 }},\ k=0,1,2,....
\end{equation}

The coefficients of the Euler-Frobenius polynomial $E_k(x)=\sum\limits_{s=0}^k\mathrm{a}_sx^s$ of degree $k$ are expressed by the following formula
which was obtained by Euler:
$$
\mathrm{a}_s=\sum\limits_{j=0}^s(-1)^j{{k+2}\choose j}(s+1-j)^{k+1}.
$$
In \cite{Frob1910} it was shown that all roots $q_j$, $j=1,2,...,k,$ of the polynomial $E_k(x)$ are real, negative and distinct, that is:
$$
q_1<q_2<...<q_k<0.
$$
Furthermore, these roots satisfy the relation
$$
q_j\cdot q_{k+1-j}=1,\ j=1,2,...,k.
$$

For the Euler-Frobenius polynomials $E_k(x)$ the following
identity holds
\begin{equation}\label{eq2.6}
E_k (x) = x^k E_k \left( {\frac{1}{x}} \right),
\end{equation}
and also the following is true.

\begin{theorem}\label{Thm1} \emph{ (Lemma 3 of \cite{Shad10})}. {Polynomial $Q_k (x)$ which is
defined by the formula}
\begin{equation}\label{eq2.7}
Q_k (x) = (x - 1)^{k + 1} \sum\limits_{i = 0}^{k + 1}
{\frac{{\Delta ^i 0^{k + 1} }}{{(x - 1)^i }}}
\end{equation}
{is the Euler-Frobenius polynomial (\ref{eq2.5}) of degree $k$, i.e.,
$Q_k(x) = E_k(x)$, where $\Delta^i0^k=\sum_{l=1}^i(-1)^{i-l}{i\choose l}l^k$.}
\end{theorem}

We also use the formula
\begin{equation}\label{eq2.8}
\sum\limits_{\gamma  = 0}^{n - 1} {q^\gamma  \gamma ^k  =
\frac{1}{{1 - q}}\sum\limits_{i = 0}^k {\left( {\frac{q}{{1 - q}}}
\right)^i \Delta ^i 0^k  - \frac{{q^n }}{{1 - q}}\sum\limits_{i =
0}^k {\left( {\frac{q}{{1 - q}}} \right)^i \Delta ^i \gamma ^k
|_{\gamma  = n} ,} } }
\end{equation}
which is given in \cite{Ham62}, where $\Delta^i\gamma^k$ is the finite difference of order $i$ of
$\gamma^k$ and $q$ is the ratio of a geometric progression.

In finding the analytic formulas for coefficients of optimal formulas in the space $L_2^{(m)}$ by Sobolev method the discrete analogue $D_m(h\beta)$ of
the operator $\frac{\d^{2m}}{\d x^{2m}}$ plays the main role. This discrete analogue  satisfies the
 equality
\begin{equation}\label{eq2.10}
hD_m(h\beta)*G_m(h\beta) = \delta_{\d}(h\beta),
\end{equation}
where $G_m(h\beta)$ is the discrete argument function for the function
\begin{equation}\label{eq2.11}
G_m(x)=\frac{|x|^{2m-1}}{2(2m-1)!},
\end{equation}
and $\delta_{\d}(h\beta)$ is equal to 0 when $\beta\neq 0$ and 1 when $\beta=0$.

We note that the operator $D_m(h\beta)$ was
introduced and studied by Sobolev \cite{Sobolev74}.
In \cite{Shad85} the discrete function $D_m(h\beta)$ was constructed and the following was
proved.

\begin{theorem}\label{Thm2}
{ The discrete analogue of the differential operator
$\frac{\d^{2m}}{\d x^{2m}}$ has the form
\begin{equation}\label{eq2.12}
D_m(h\beta)=p\left\{
\begin{array}{lll}
{\displaystyle \sum\limits_{k=1}^{m-1}A_kq_k^{|\beta|-1}        }
& \mbox{ for }& |\beta|\geq 2,\\
{\displaystyle 1+\sum\limits_{k=1}^{m-1}A_k                  }
& \mbox{ for }& |\beta|= 1,\\
{\displaystyle C+\sum\limits_{k=1}^{m-1}\frac{A_k}{q_k} } & \mbox{
for }& \beta= 0,
\end{array}
\right.
\end{equation}
where
\begin{equation}\label{eq2.13}
p=\frac{(2m-1)!}{h^{2m}},\quad
A_k=\frac{(1-q_k)^{2m+1}}{E_{2m-1}(q_k)}, \quad C=-2^{2m-1},
\end{equation}
$E_{2m-1}(x)$ is the Euler-Frobenius polynomial of degree $2m-1$,
$q_k$ are the roots of the Euler-Frobenius polynomial
$E_{2m-2}(x)$, $|q_k|<1$, and $h$ is a small positive parameter.}
\end{theorem}

In addition, some properties of $D_m(h\beta)$ were studied in the works \cite{Shad85,Sobolev74}. Here we give the
following.

\begin{theorem}\label{Thm3} { The discrete argument function $D_m(h\beta)$ and the
monomials $(h\beta)^k$ are related to each other as follows}:
\begin{equation}\label{eq2.14}
\sum_{\beta=-\infty}^{\infty}D_m(h\beta)(h\beta)^k= \left\{
\begin{array}{lll}
0&\mbox{\it for } & 0\leq k\leq 2m-1,\\
(2m)!&\mbox{\it for } & k= 2m.
\end{array} \right.
\end{equation}
\end{theorem}

Now we give the results of the works \cite{Shad83,Shad10} on the optimal quadrature formulas
of the form (\ref{eq2.1}) in the sense of Sard and on the norm of the optimal error functional corresponding to the case $\omega=0$.

\begin{theorem}\label{Thm4} Coefficients of the optimal quadrature formulas of the form
\begin{equation}\label{eq2.15}
\int\limits_0^1\varphi(x)\d x\cong \sum\limits_{\beta=0}^NC_{\beta}\varphi(h\beta)
\end{equation}
in the space $L_2^{(m)}[0,1]$ have the form
\begin{eqnarray}
\mathring{C}_0&=&h\left(\frac{1}{2}-\sum\limits_{k=1}^{m-1}d_k\frac{q_k-q_k^N}{1-q_k}\right),\nonumber\\
\mathring{C}_{\beta}&=&h\left(1+\sum\limits_{k=1}^{m-1}d_k\left(q_k^{\beta}+q_k^{N-\beta}\right)\right),\ \beta=1,2,...,N-1,\label{eq2.16}\\
\mathring{C}_N&=&h\left(\frac{1}{2}-\sum\limits_{k=1}^{m-1}d_k\frac{q_k-q_k^N}{1-q_k}\right),\nonumber
\end{eqnarray}
where $d_k$ satisfy the system
$$
\sum\limits_{k=1}^{m-1}d_k\sum\limits_{i=1}^j\frac{q_k+(-1)^{i+1}q_k^{N+i}}{(q_k-1)^{i+1}}\Delta^i0^j=\frac{B_{j+1}}{j+1},\ j=1,2,...,m-1,
$$
$q_k$ are roots of the Euler-Frobenius polynomial $E_{2m-2}(x)$ of degree $2m-2$ with $|q_k|<1$, and $\Delta^i0^j=\sum_{l=1}^i(-1)^{i-l}{i\choose l}l^j$.
\end{theorem}

\begin{theorem}\label{Thm5}
The norm of the error functional of the optimal quadrature formula (\ref{eq2.15}) on the space $L_2^{(m)}[0,1]$
has the form
\begin{equation*}
\|\mathring{\ell}\|^2_{L_2^{(m)*}[0,1]}=(-1)^{m+1}
\left(
\frac{h^{2m}B_{2m}}{(2m)!}+\frac{2h^{2m+1}}{(2m)!}\sum\limits_{k=1}^{m-1}d_k\sum\limits_{i=1}^{2m}\frac{-q_k^{N+i}+(-1)^iq_k}{(1-q_k)^{i+1}}\Delta^i0^{2m}
\right),
\end{equation*}
where $B_{2m}$ is the Bernoulli number, and $q_k$, $d_k$ and $\Delta^i0^{2m}$ are given in Theorem~\ref{Thm4}.
\end{theorem}

\subsection{Construction of optimal quadrature formulas for the interval $[0,1]$}

Here we obtain optimal quadrature formulas of the form (\ref{eq2.1}) for the interval $[0,1]$ when $\omega\in \mathbb{R}$
and $\omega\neq 0$. In the space $L_2^{(m)}[0,1]$, using the results of Sections 2, 3 and 5 of \cite{BolHayShad17},
for the coefficients of the optimal quadrature formulas in the sense of Sard of the form
\begin{equation}\label{eq2.17}
\int\limits_0^1 {\E^{2\pi \i\omega x}\varphi (x)\d x \cong } \sum\limits_{\beta  = 0}^N
{C_\beta  } \varphi (h\beta  )
\end{equation}
for $\omega\in \mathbb{R}$ with $\omega\neq 0$, we get the following system of linear equations
\begin{eqnarray}
&&\sum\limits_{\gamma=0}^N C_\gamma G_m(h\beta   - h\gamma) +P_{m-1}(h\beta) =f_m(h\beta), \ \beta  =0,1,...,N, \label{eq2.18}\\
&&\sum\limits_{\gamma=0}^N C_\gamma
(h\gamma)^{\alpha}=g_{\alpha},\ \alpha=0,1,...,m-1,\label{eq2.19}
\end{eqnarray}
where $P_{m-1}(h\beta)=\sum\limits_{\alpha=0}^{m-1}p_{\alpha}(h\beta)^{\alpha}$ is a polynomial of degree $m-1$
with complex coefficients,
\begin{eqnarray}
f_m(h\beta)&=& \int\limits_0^1\E^{2\pi \i\omega x} G_m(x - h\beta  )\d x \label{eq2.20}\\
&=&-\sum\limits_{\alpha=0}^{2m-1}\frac{(h\beta)^{2m-1-\alpha}(-1)^{\alpha}g_{\alpha}}{2\alpha!(2m-1-\alpha)!}
+\frac{\E^{2\pi\i\omega h\beta}}{(2\pi\i\omega)^{2m}}-\sum\limits_{k=0}^{2m-1}\frac{(h\beta)^{2m-1-k}}{(2m-1-k)!(2\pi\i\omega)^{k+1}}, \nonumber \\
g_{\alpha}&=&\int\limits_0^1\E^{2\pi \i \omega x}x^{\alpha}\d x \label{eq2.21} \\
&=&\sum\limits_{k=0}^{\alpha-1}\frac{(-1)^k\alpha!\ \E^{2\pi\i\omega}}{(\alpha-k)!(2\pi\i\omega)^{k+1}}+\frac{(-1)^{\alpha}\alpha!}{(2\pi \i \omega)^{\alpha+1}}\left(\E^{2\pi\i\omega}-1\right),\ \alpha=0,1,..., \nonumber
\end{eqnarray}
$G_m(x)$ is defined by (\ref{eq2.11}), and $h=\frac{1}{N}$ for $N\in \mathbb{N}$.

In the system (\ref{eq2.18})-(\ref{eq2.19}) unknowns are the optimal coefficients $\mathring{C}_\beta,$ $\beta=0,1,...,N$ and $p_{\alpha}$, $\alpha=0,1,...,m-1$. We point out that when $N+1\geq m$ the system (\ref{eq2.18})-(\ref{eq2.19}) has a unique solution.
This solution satisfies conditions (\ref{eq2.3-1}) and the equality (\ref{eq2.3-2}). It should be noted that
the existence and uniqueness of the solution for such type of systems were studied, for example,
in \cite{ShadHayAkhm15,Sobolev74,SobVas}.

We are interested in finding explicit formulas for the optimal coefficients $\mathring{C}_{\beta}$, $\beta=0,1,...,N$ and unknown polynomial $P_{m-1}(h\beta)$ satisfying the system (\ref{eq2.18})-(\ref{eq2.19}). The system (\ref{eq2.18})-(\ref{eq2.19}) is solved similarly as the system (34)-(35) of \cite{BolHayShad17} by Sobolev's method, using the discrete analogue $D_m(h\beta)$ of the differential operator $\frac{\d^{2m}}{\d x^{2m}}$.

We formulate the results of this section as the following two theorems.

\begin{theorem}\label{Thm6}
{For real $\omega$ with $\omega h\not\in \mathbb{Z}$, the coefficients of optimal quadrature
formulas of the form (\ref{eq2.1})  in the space
$L_2^{(m)}[0,1]$ when $N+1\geq m$ are expressed by formulas}
\begin{eqnarray}
\mathring{C}_0&=&h\left(\frac{\E^{2\pi \i\omega h}K_{\omega,m}}{\E^{2\pi \i\omega h}-1}-\frac{1}{2\pi \i\omega h}+\sum\limits_{k=1}^{m-1}\left(a_k\frac{q_k}{q_k-1}+b_k\frac{q_k^N}{1-q_k}\right)\right), \label{eq2.22}  \\
\mathring{C}_\beta&=&h\left(\E^{2\pi \i\omega h\beta}K_{\omega,m}+\sum\limits_{k=1}^{m-1}\left(a_kq_k^\beta+b_kq_k^{N-\beta}\right)\right),
\ \beta=1,2,...,N-1, \label{eq2.23} \\
\mathring{C}_N&=&h\left(\frac{\E^{2\pi \i\omega}K_{\omega,m}}{1-\E^{2\pi \i\omega h}}+\frac{\E^{2\pi \i\omega}}{2\pi \i\omega h}+\sum\limits_{k=1}^{m-1}\left(a_k\frac{q_k^N}{1-q_k}+b_k\frac{q_k}{q_k-1}\right)\right),  \label{eq2.24}
\end{eqnarray}
{where $a_k$ and $b_k$ are defined by the following system of $(2m-2)$ linear equations}\\
\begin{equation}  \label{eq2.25}
\begin{array}{l}
\sum\limits_{k=1}^{m-1}a_k\left[\sum\limits_{t=1}^{j}\frac{q_k\Delta^t0^j}{(q_k-1)^{t+1}}\right]+
\sum\limits_{k=1}^{m-1}b_k\left[\sum\limits_{t=1}^j\frac{q_k^{N+t}\Delta^t0^{j}}{(1-q_k)^{t+1}}\right]\\
\quad =\frac{j!}{(2\pi \i\omega h)^{j+1}}-\sum\limits_{t=1}^j\frac{\E^{2\pi \i\omega h}K_{\omega,m}\Delta^t0^j}{(\E^{2\pi \i\omega h}-1)^{t+1}},                 \ j=1,2,...,m-1,\\
\sum\limits_{k=1}^{m-1}a_k\Bigg[\sum\limits_{t=1}^j\frac{q_k^t\Delta^t0^j}{(1-q_k)^{t+1}}-\sum\limits_{\alpha=1}^{j}h^{\alpha-j}{j\choose {\alpha}}
\sum\limits_{t=1}^{\alpha}\frac{q_k^{N+t}\Delta^t0^{\alpha}}{(1-q_k)^{t+1}}\Bigg]\\
\quad+\sum\limits_{k=1}^{m-1}b_k\Bigg[\sum\limits_{t=1}^j\frac{q_k^{N+1}\Delta^t0^j}{(q_k-1)^{t+1}}-\sum\limits_{\alpha=1}^{j}h^{\alpha-j}{j\choose {\alpha}}
\sum\limits_{t=1}^{\alpha}\frac{q_k\Delta^t0^{\alpha}}{(q_k-1)^{t+1}}\Bigg]\\
\quad=\frac{(-1)^{j+1}j!}{(2\pi \i\omega h)^{j+1}}
+\sum\limits_{\alpha=1}^{j}h^{\alpha-j}\frac{(-1)^{\alpha}j!\E^{2\pi \i\omega}}{(j-\alpha)!(2\pi \i\omega h)^{\alpha+1}}
-\frac{K_{\omega,m}}{1-\E^{2\pi \i\omega h}}\sum\limits_{t=1}^j
\left(\frac{\E^{2\pi \i\omega h}}{1-\E^{2\pi \i\omega h}}\right)^t\Delta^t0^j\\
\quad\quad+\frac{K_{\omega,m}\E^{2\pi \i\omega}}{1-\E^{2\pi \i\omega h}}\sum\limits_{\alpha=1}^{j}h^{\alpha-j}{j\choose {\alpha}}\sum\limits_{t=1}^{\alpha}
\left(\frac{\E^{2\pi \i\omega h}}{1-\E^{2\pi \i\omega h}}\right)^t\Delta^t0^\alpha,\
\ j=1,2,...,m-1,
\end{array}
\end{equation}
{$q_k$ are roots of the Euler-Frobenius polynomial $E_{2m-2}(x)$
with $|q_k|<1$,} and
\begin{equation} \label{eq2.26}
K_{\omega,m}=\left(\frac{\sin \pi\omega h}{\pi\omega h}\right)^{2m}\frac{(2m-1)!}{2\sum\limits_{\alpha=0}^{m-2}\mathrm{a}_{\alpha}\cos[2\pi\omega h(m-1-\alpha)]+\mathrm{a}_{m-1}}.
\end{equation}
Here $\mathrm{a}_{\alpha}=\sum\limits_{j=0}^{\alpha}(-1)^j{{2m}\choose j}(s+1-j)^{2m-1}$ are the coefficients of the Euler-Frobenius polynomial $E_{2m-2}(x)$
of degree $2m-2$.
\end{theorem}

\begin{theorem}\label{Thm7} {For $\omega h\in \mathbb{Z}$ with $\omega\neq 0$, the coefficients of optimal quadrature
formulas of the form (\ref{eq2.1}) in the space
$L_2^{(m)}[0,1]$ when $N+1\geq m$ are expressed by formulas}
\begin{eqnarray*}
\mathring{C}_0&=&h\left(-\frac{1}{2\pi \i\omega h}+\sum\limits_{k=1}^{m-1}\left(a_k\frac{q_k}{q_k-1}+b_k\frac{q_k^N}{1-q_k}\right)\right),\\
\mathring{C}_\beta&=&h\sum\limits_{k=1}^{m-1}\left(a_kq_k^\beta+b_kq_k^{N-\beta}\right),
\ \ \ \ \beta=1,2,...,N-1,\\
\mathring{C}_N&=&h\left(\frac{\E^{2\pi\i\omega}}{2\pi \i\omega h}+\sum\limits_{k=1}^{m-1}\left(a_k\frac{q_k^N}{1-q_k}+b_k\frac{q_k}{q_k-1}\right)\right),
\end{eqnarray*}
{where $a_k$ and $b_k$, $k=1,2,...,m-1$, are defined by the following system of $(2m-2)$ linear equations:}\\
$$
\begin{array}{l}
\sum\limits_{k=1}^{m-1}a_k\left[\sum\limits_{t=1}^{j}\frac{q_k\Delta^t0^j}{(q_k-1)^{t+1}}\right]+
\sum\limits_{k=1}^{m-1}b_k\left[\sum\limits_{t=1}^j\frac{q_k^{N+t}\Delta^t0^{j}}{(1-q_k)^{t+1}}\right]=\frac{j!}{(2\pi \i\omega h)^{j+1}},\ j=1,2,...,m-1,\\
\sum\limits_{k=1}^{m-1}a_k\Bigg[\sum\limits_{t=1}^j\frac{q_k^t\Delta^t0^j}{(1-q_k)^{t+1}}-\sum\limits_{\alpha=1}^{j}h^{\alpha-j}{j\choose {\alpha}}
\sum\limits_{t=1}^{\alpha}\frac{q_k^{N+t}\Delta^t0^{\alpha}}{(1-q_k)^{t+1}}\Bigg]\\
+\sum\limits_{k=1}^{m-1}b_k\Bigg[\sum\limits_{t=1}^j\frac{q_k^{N+1}\Delta^t0^j}{(q_k-1)^{t+1}}
-\sum\limits_{\alpha=1}^{j}h^{\alpha-j}{j\choose {\alpha}}
\sum\limits_{t=1}^{\alpha}\frac{q_k\Delta^t0^{\alpha}}{(q_k-1)^{t+1}}\Bigg]\\
\qquad=\frac{(-1)^{j+1}j!}{(2\pi \i\omega h)^{j+1}}
+\sum\limits_{\alpha=1}^{j}h^{\alpha-j}\frac{(-1)^{\alpha}j!\E^{2\pi \i\omega}}{(j-\alpha)!(2\pi \i\omega h)^{\alpha+1}},\
\ j=1,2,...,m-1.
\end{array}
$$
{Here $q_k$ are the roots of the Euler-Frobenius polynomial $E_{2m-2}(x)$
 and $|q_k|<1$.}
\end{theorem}

We note that Theorem \ref{Thm6} is generalization of Theorem 6 in [6] for real $\omega$ with $\omega h\not\in \mathbb{Z}$
while Theorem \ref{Thm7} for $\omega h\in \mathbb{Z}$ with $\omega\neq 0$ is the same with Theorem 7 of the work [6]. Theorem 6 is proved similarly as Theorem 6 of [6].
Therefore, it is sufficient to give a brief proof of Theorem \ref{Thm6}.

\emph{The brief proof of Theorem \ref{Thm6}}. First, such as in Theorem 5 of the work \cite{BolHayShad17}, using the discrete function $D_m(h\beta)$, for
optimal coefficients $\mathring{C}_{\beta},$ $\beta=1,2,...,N-1$, when $\omega$ is real and $\omega h\not\in \mathbb{Z}$ we get the following formula
\begin{equation}\label{eq2.27}
\mathring{C}_\beta=h\left(K_{\omega,m}\E^{2\pi\i \omega h\beta}+\sum\limits_{k=1}^{m-1}
\left(a_kq_k^{\beta}+b_kq_k^{N-\beta}\right)\right) ,\ \
\beta=1,2,...,N-1,
\end{equation}
where $a_k$, $b_k$, $k=1,2,...,m-1$, $K_{\omega,m}$ are unknowns,
and $q_k$ are roots of the Euler-Frobenius polynomial $E_{2m-2}(x)$ of degree $2m-2$ with $|q_k|<1$.

Next, it is sufficient to find $a_k,\ b_k$, $k=1,2,...,m-1$, $K_{\omega,m}$,
optimal coefficients $\mathring{C}_0$, $\mathring{C}_N$ and unknown polynomial $P_{m-1}(h\beta)$ of degree $m-1$.
Now putting the form (\ref{eq2.27}) of optimal coefficients $\mathring{C}_{\beta}$, $\beta=1,2,...,N-1$, into (\ref{eq2.18}), using (\ref{eq2.6})-(\ref{eq2.8}), after some simplifications,
we get the following identity with respect to $(h\beta)$
$$
\begin{array}{l}
\E^{2\pi\i\omega h\beta}\left(\frac{h^{2m}K_{\omega,m}}{(2m-1)!}\sum\limits_{i=0}^{2m-1}\frac{\E^{2\pi\i\omega h}\Delta^i0^{2m-1}}{(\E^{2\pi\i\omega h}-1)^{i+1}}\right)\\
+\frac{(h\beta)^{2m-1}}{(2m-1)!}\left(\mathring{C}_0-\frac{K_{\omega,m}h\E^{2\pi\i\omega h}}{\E^{2\pi\i\omega h}-1}-h\sum\limits_{k=1}^{m-1}\frac{a_kq_k-b_kq_k^N}{q_k-1}\right)\\
-\sum\limits_{j=1}^{2m-1}\frac{(h\beta)^{2m-1-j}}{j!(2m-1-j)!}h^{j+1}\left(\sum\limits_{i=0}^j\frac{K_{\omega,m}\E^{2\pi\i\omega h}\Delta^i0^j}
{(\E^{2\pi\i\omega h}-1)^{i+1}}+\sum\limits_{k=1}^{m-1}\sum\limits_{i=0}^j\frac{a_kq_k+(-1)^{i+1}b_kq_k^{N+i}}{(q_k-1)^{i+1}}\Delta^i0^j\right)\\
-\sum\limits_{j=m}^{2m-1}\frac{(h\beta)^{2m-1-j}}{j!(2m-1-j)!}\frac{(-1)^j}{2}\sum\limits_{\gamma=0}^{N}\mathring{C}_{\gamma}(h\gamma)^j+P_{m-1}(h\beta)\\
\qquad\qquad=\frac{\E^{2\pi\i\omega h\beta}}{(2\pi\i\omega)^{2m}}-\frac{(h\beta)^{2m-1}}{(2m-1)!(2\pi \i \omega)}\\
\qquad\qquad\quad-\sum\limits_{j=1}^{2m-1}\frac{(h\beta)^{2m-1-j}}{j!(2m-1-j)!}\frac{j!}{(2\pi\i\omega)^{j+1}}
-\sum\limits_{j=m}^{2m-1}\frac{(h\beta)^{2m-1-j}}{j!(2m-1-j)!}\frac{(-1)^j}{2}g_j.
\end{array}
$$
Hence, the formulas (\ref{eq2.22}) and (\ref{eq2.26}) for $\mathring{C}_0$ and $K_{\omega,m}$ are obtained by equating the coefficients of the terms $(h\beta)^{2m-1}$ and  $\E^{2\pi\i\omega h\beta}$, respectively. Then (\ref{eq2.24}) is found from (\ref{eq2.19}) when $\alpha=0$, using (\ref{eq2.22}) and
(\ref{eq2.27}). From the last identity equating the coefficients of the terms $(h\beta)^{2m-1-j}$ for $j=1,2,...,m-1$ and using (\ref{eq2.19}) for
$\alpha=1,2,...,m-1$, taking into account (\ref{eq2.22}),  (\ref{eq2.24}) and (\ref{eq2.27}), we obtain the system
(\ref{eq2.25}). Finally, the polynomial $P_{m-1}(h\beta)$ is obtained from the identity by equating  the coefficients of the terms
$(h\beta)^{2m-1-j}$ for $j=m,m+1,...,2m-1$. Theorem is proved. \hfill $\Box$


\subsection{Optimal quadrature formulas for the interval $[a,b]$}

In this section we obtain the optimal quadrature formulas for the interval $[a,b]$ in the space $L_2^{(m)}$
by a linear transformation from the results of the previous section.

We consider construction of the optimal quadrature formulas of the form
\begin{equation}\label{eq2.28}
\int\limits_a^b\E^{2\pi \i \omega x}\varphi(x)\ \d x\cong \sum\limits_{\beta=0}^NC_{\beta,\omega}[a,b]\varphi(x_\beta)
\end{equation}
in the Sobolev space
$L_2^{(m)}[a,b]$. Here $C_{\beta,\omega}[a,b]$ are coefficients,
$x_\beta=h\beta+a$ $(\in [a,b])$ are nodes of the formula (\ref{eq2.28}),
$\omega\in \mathbb{R}$, $\i^2=-1$, and $h=\frac{b-a}{N}$ for $N+1\geq m$.

Now, by a linear transformation $x=(b-a)y+a$, where $0\leq y\leq 1$, we obtain
\begin{equation}\label{eq2.29}
\int\limits_a^b\E^{2\pi \i \omega x}\varphi(x)\ \d x=(b-a)\E^{2\pi \i \omega a}\int\limits_0^1\E^{2\pi \i \omega (b-a)y}\varphi((b-a)y+a)\d y.
\end{equation}
Then, applying Theorems \ref{Thm4}, \ref{Thm6} and \ref{Thm7} to the integral on the right-hand side of the last equality,
we have the following results which are optimal quadrature formulas of the form (\ref{eq2.28}) in the sense of Sard
in the space $L_2^{(m)}[a,b]$ for all real $\omega$.

For the case $\omega=0$, using Theorem \ref{Thm4} in (\ref{eq2.29}), we get

\begin{theorem}\label{Thm8} Coefficients of the optimal quadrature formulas of the form
\begin{equation*}
\int\limits_a^b\varphi(x)\d x\cong \sum\limits_{\beta=0}^NC_{\beta,0}[a,b]\varphi(h\beta+a)
\end{equation*}
in the space $L_2^{(m)}[a,b]$ have the form
\begin{eqnarray*}
\mathring{C}_{0,0}[a,b]&=&h\left(\frac{1}{2}-\sum\limits_{k=1}^{m-1}d_k\frac{q_k-q_k^N}{1-q_k}\right),\nonumber\\
\mathring{C}_{\beta,0}[a,b]&=&h\left(1+\sum\limits_{k=1}^{m-1}d_k\left(q_k^{\beta}+q_k^{N-\beta}\right)\right),\ \beta=1,2,...,N-1,\label{eq2.16}\\
\mathring{C}_{N,0}[a,b]&=&h\left(\frac{1}{2}-\sum\limits_{k=1}^{m-1}d_k\frac{q_k-q_k^N}{1-q_k}\right),\nonumber
\end{eqnarray*}
where $h=\frac{b-a}{N}$ and $d_k$ satisfy the system
$$
\sum\limits_{k=1}^{m-1}d_k\sum\limits_{i=1}^j\frac{q_k+(-1)^{i+1}q_k^{N+i}}{(q_k-1)^{i+1}}\Delta^i0^j=\frac{B_{j+1}}{j+1},\ j=1,2,...,m-1.
$$
Here $B_{j+1}$ is the Bernoulli number, $q_k$ are roots of the Euler-Frobenius polynomial $E_{2m-2}(x)$ of degree $(2m-2)$ with $|q_k|<1$, and $\Delta^i0^j=\sum_{l=1}^i(-1)^{i-l}{i\choose l}l^j$.
\end{theorem}

For the case $\omega\in \mathbb{R}$ and $\omega h\not\in \mathbb{Z}$, applying Theorem \ref{Thm6} to the integral in (\ref{eq2.29}), we obtain

\begin{theorem}\label{Thm9}
{For real $\omega$ with $\omega h\not\in \mathbb{Z}$, the coefficients of optimal quadrature
formulas of the form (\ref{eq2.28})  in the space
$L_2^{(m)}[a,b]$ when $N+1\geq m$ are expressed by formulas}
\begin{eqnarray*}
\mathring{C}_{0,\omega}[a,b]&=&h\left(\frac{\E^{2\pi \i\omega (a+h)}K_{\omega,m}}{\E^{2\pi \i\omega h}-1}-\frac{\E^{2\pi \i\omega a}}{2\pi \i\omega h}+\sum\limits_{k=1}^{m-1}\left(a_k\frac{q_k}{q_k-1}+b_k\frac{q_k^N}{1-q_k}\right)\right),   \\
\mathring{C}_{\beta,\omega}[a,b]&=&h\left(\E^{2\pi \i\omega (h\beta+a)}K_{\omega,m}+\sum\limits_{k=1}^{m-1}\left(a_kq_k^\beta+b_kq_k^{N-\beta}\right)\right),
\ \beta=1,2,...,N-1,  \\
\mathring{C}_{N,\omega}[a,b]&=&h\left(\frac{\E^{2\pi \i\omega b}K_{\omega,m}}{1-\E^{2\pi \i\omega h}}+\frac{\E^{2\pi \i\omega b}}{2\pi \i\omega h}+\sum\limits_{k=1}^{m-1}\left(a_k\frac{q_k^N}{1-q_k}+b_k\frac{q_k}{q_k-1}\right)\right),
\end{eqnarray*}
{where $a_k$ and $b_k$ are defined by the following system of $(2m-2)$ linear equations:}\\
\begin{equation*}
\begin{array}{l}
\sum\limits_{k=1}^{m-1}a_k\left[\sum\limits_{t=1}^{j}\frac{q_k\Delta^t0^j}{(q_k-1)^{t+1}}\right]+
\sum\limits_{k=1}^{m-1}b_k\left[\sum\limits_{t=1}^j\frac{q_k^{N+t}\Delta^t0^{j}}{(1-q_k)^{t+1}}\right]\\
\qquad =\frac{j!\E^{2\pi \i\omega a}}{(2\pi \i\omega h)^{j+1}}-\sum\limits_{t=1}^j\frac{\E^{2\pi \i\omega (a+h)}K_{\omega,m}\Delta^t0^j}{(\E^{2\pi \i\omega h}-1)^{t+1}},                 \ j=1,2,...,m-1,\\
\sum\limits_{k=1}^{m-1}a_k\Bigg[\sum\limits_{t=1}^j\frac{q_k^t\Delta^t0^j}{(1-q_k)^{t+1}}-\sum\limits_{\alpha=1}^{j}(\frac{h}{b-a})^{\alpha-j}{j\choose {\alpha}}
\sum\limits_{t=1}^{\alpha}\frac{q_k^{N+t}\Delta^t0^{\alpha}}{(1-q_k)^{t+1}}\Bigg]\\
\quad+\sum\limits_{k=1}^{m-1}b_k\Bigg[\sum\limits_{t=1}^j\frac{q_k^{N+1}\Delta^t0^j}{(q_k-1)^{t+1}}-\sum\limits_{\alpha=1}^{j}(\frac{h}{b-a})^{\alpha-j}{j\choose {\alpha}}
\sum\limits_{t=1}^{\alpha}\frac{q_k\Delta^t0^{\alpha}}{(q_k-1)^{t+1}}\Bigg]\\
\quad=\frac{(-1)^{j+1}j!\E^{2\pi \i\omega a}}{(2\pi \i\omega h)^{j+1}}
+\sum\limits_{\alpha=1}^{j}(\frac{h}{b-a})^{\alpha-j}\frac{(-1)^{\alpha}j!\E^{2\pi \i\omega b}}{(j-\alpha)!(2\pi \i\omega h)^{\alpha+1}}
-\frac{K_{\omega,m}\E^{2\pi \i\omega a}}{1-\E^{2\pi \i\omega h}}\sum\limits_{t=1}^j
\left(\frac{\E^{2\pi \i\omega h}}{1-\E^{2\pi \i\omega h}}\right)^t\Delta^t0^j\\
\quad\quad+\frac{K_{\omega,m}\E^{2\pi \i\omega b}}{1-\E^{2\pi \i\omega h}}\sum\limits_{\alpha=1}^{j}(\frac{h}{b-a})^{\alpha-j}{j\choose {\alpha}}\sum\limits_{t=1}^{\alpha}
\left(\frac{\E^{2\pi \i\omega h}}{1-\E^{2\pi \i\omega h}}\right)^t\Delta^t0^\alpha,\
\ j=1,2,...,m-1,
\end{array}
\end{equation*}
{$q_k$ are roots of the Euler-Frobenius polynomial $E_{2m-2}(x)$
with $|q_k|<1$,} and $K_{\omega,m}$ is defined by (\ref{eq2.26}).
\end{theorem}

Lastly, for the case $\omega h\in \mathbb{Z}$ with $\omega\neq 0$ application Theorem \ref{Thm7} in (\ref{eq2.29}) gives the following.

\begin{theorem}\label{Thm10} {For $\omega h\in \mathbb{Z}$ with $\omega\neq 0$, the coefficients of optimal quadrature
formulas of the form (\ref{eq2.28}) in the space
$L_2^{(m)}[a,b]$ when $N+1\geq m$ are expressed by formulas}
\begin{eqnarray*}
\mathring{C}_{0,\omega}[a,b]&=&h\left(-\frac{\E^{2\pi\i\omega a}}{2\pi \i\omega h}+\sum\limits_{k=1}^{m-1}\left(a_k\frac{q_k}{q_k-1}+b_k\frac{q_k^N}{1-q_k}\right)\right),\\
\mathring{C}_{\beta,\omega}[a,b]&=&h\sum\limits_{k=1}^{m-1}\left(a_kq_k^\beta+b_kq_k^{N-\beta}\right),
\ \ \ \ \beta=1,2,...,N-1,\\
\mathring{C}_{N,\omega}[a,b]&=&h\left(\frac{\E^{2\pi\i\omega b}}{2\pi \i\omega h}+\sum\limits_{k=1}^{m-1}\left(a_k\frac{q_k^N}{1-q_k}+b_k\frac{q_k}{q_k-1}\right)\right),
\end{eqnarray*}
{where $a_k$ and $b_k$, $k=1,2,...,m-1$, are defined by the following system of $(2m-2)$ linear equations}\\
$$
\begin{array}{l}
\sum\limits_{k=1}^{m-1}a_k\left[\sum\limits_{t=1}^{j}\frac{q_k\Delta^t0^j}{(q_k-1)^{t+1}}\right]+
\sum\limits_{k=1}^{m-1}b_k\left[\sum\limits_{t=1}^j\frac{q_k^{N+t}\Delta^t0^{j}}{(1-q_k)^{t+1}}\right]=\frac{j!\E^{2\pi \i\omega a}}{(2\pi \i\omega h)^{j+1}},\ j=1,2,...,m-1,\\
\sum\limits_{k=1}^{m-1}a_k\Bigg[\sum\limits_{t=1}^j\frac{q_k^t\Delta^t0^j}{(1-q_k)^{t+1}}-\sum\limits_{\alpha=1}^{j}(\frac{h}{b-a})^{\alpha-j}{j\choose {\alpha}}
\sum\limits_{t=1}^{\alpha}\frac{q_k^{N+t}\Delta^t0^{\alpha}}{(1-q_k)^{t+1}}\Bigg]\\
\qquad+\sum\limits_{k=1}^{m-1}b_k\Bigg[\sum\limits_{t=1}^j\frac{q_k^{N+1}\Delta^t0^j}{(q_k-1)^{t+1}}
-\sum\limits_{\alpha=1}^{j}(\frac{h}{b-a})^{\alpha-j}{j\choose {\alpha}}
\sum\limits_{t=1}^{\alpha}\frac{q_k\Delta^t0^{\alpha}}{(q_k-1)^{t+1}}\Bigg]\\
\qquad\qquad=
\frac{(-1)^{j+1}j!\E^{2\pi \i\omega a}}{(2\pi \i\omega h)^{j+1}}
+\sum\limits_{\alpha=1}^{j}(\frac{h}{b-a})^{\alpha-j}\frac{(-1)^{\alpha}j!\E^{2\pi \i\omega b}}{(j-\alpha)!(2\pi \i\omega h)^{\alpha+1}},\
\ j=1,2,...,m-1.
\end{array}
$$
{Here $q_k$ are the roots of the Euler-Frobenius polynomial $E_{2m-2}(x)$
 and $|q_k|<1$.}
\end{theorem}

From a practical point of view, to have the formulas for coefficients of the optimal quadrature formulas
(\ref{eq2.28}) for the first several values of $m$ is very useful.
Below we present explicit formulas of the optimal coefficients for $m=2$ and 3, using Theorems \ref{Thm8}, \ref{Thm9} and \ref{Thm10}.

\textbf{Remark 1.}  In the space $L_2^{(1)}[a,b]$ the coefficients of the optimal quadrature formulas (\ref{eq2.28})
have the same forms as derived in \cite{HayJeonLee19}, where it was shown that for functions with a continuous second derivative
the convergence order of this optimal quadrature formula is $O(h^2)$.

\textbf{Remark 2.} The optimal quadrature formulas (\ref{eq2.28}) for $m=1$ and 2 have the same convergence order.

\begin{corollary}\label{Cor2} In the space $L_2^{(2)}[a,b]$ the coefficients of the optimal quadrature formulas (\ref{eq2.28})
for $\omega=0$ are written as
\begin{eqnarray*}
\mathring{C}_{0,0}[a,b]&=&h\left(\frac{1}{2}+\frac{q_1-q_1^N}{2(1-q_1)(1+q_1^N)}\right),\\
\mathring{C}_{\beta,0}[a,b]&=&h\left(1-\frac{q_1^{\beta}+q_1^{N-\beta}}{2(1+q_1^N)}\right),\ \beta=1,2,...,N-1,\\
\mathring{C}_{N,0}[a,b]&=&h\left(\frac{1}{2}+\frac{q_1-q_1^N}{2(1-q_1)(1+q_1^N)}\right),
\end{eqnarray*}
for $\omega\in \mathbb{R}$ and $\omega h \not\in \mathbb{Z}$ are expressed as
\begin{eqnarray*}
\mathring{C}_{0,\omega}[a,b]&=&h\left(\frac{\E^{2\pi \i\omega (a+h)}K_{\omega,2}}{\E^{2\pi \i\omega h}-1}-\frac{\E^{2\pi \i\omega a}}{2\pi \i\omega h}+a_1\frac{q_1}{q_1-1}+b_1\frac{q_1^N}{1-q_1}\right),   \\
\mathring{C}_{\beta,\omega}[a,b]&=&h\left(\E^{2\pi \i\omega (h\beta+a)}K_{\omega,2}+
a_1q_1^\beta+b_1q_1^{N-\beta}\right),
\ \beta=1,2,...,N-1,  \\
\mathring{C}_{N,\omega}[a,b]&=&h\left(\frac{\E^{2\pi \i\omega b}K_{\omega,2}}{1-\E^{2\pi \i\omega h}}+\frac{\E^{2\pi \i\omega b}}{2\pi \i\omega h}+a_1\frac{q_1^N}{1-q_1}+b_1\frac{q_1}{q_1-1}\right),
\end{eqnarray*}
and for $\omega h\in \mathbb{Z}$ with $\omega \neq 0$ have the form
\begin{eqnarray*}
\mathring{C}_{0,\omega}[a,b]&=&h\left(-\frac{\E^{2\pi \i\omega a}}{2\pi \i\omega h}+a_1\frac{q_1}{q_1-1}+b_1\frac{q_1^N}{1-q_1}\right),   \\
\mathring{C}_{\beta,\omega}[a,b]&=&h\left(a_1q_1^\beta+b_1q_1^{N-\beta}\right),
\ \beta=1,2,...,N-1,  \\
\mathring{C}_{N,\omega}[a,b]&=&h\left(\frac{\E^{2\pi \i\omega b}}{2\pi \i\omega h}+a_1\frac{q_1^N}{1-q_1}+b_1\frac{q_1}{q_1-1}\right),
\end{eqnarray*}
where $q_1=\sqrt{3}-2$,
$$
\begin{array}{rlrl}
a_1&= \frac{B_{\omega}(\E^{2\pi\i \omega a}-\E^{2\pi\i\omega b}q_1^N)}{1-q_1^{2N}},&
b_1&=  \frac{B_{\omega}(\E^{2\pi\i \omega b}-\E^{2\pi\i\omega a}q_1^N)}{1-q_1^{2N}},\\
B_{\omega}&=  6\left(\frac{1}{(2\pi\omega h)^2}-\frac{K_{\omega,2}}{2-2\cos 2\pi \omega h}\right),&
K_{\omega,2}&= \left(\frac{\sin \pi\omega h}{\pi\omega h}\right)^4\frac{3}{2+\cos 2\pi\omega h}.
\end{array}
$$
\end{corollary}

\begin{corollary}\label{Cor3} In the space $L_2^{(3)}[a,b]$ the coefficients of the optimal quadrature formulas (\ref{eq2.28})
for $\omega=0$ are written as
\begin{eqnarray*}
\mathring{C}_{0,0}[a,b]&=&h\left(\frac{1}{2}-\sum\limits_{k=1}^2d_k\frac{q_k-q_k^N}{1-q_k}\right),\\
\mathring{C}_{\beta,0}[a,b]&=&h\left(1+\sum\limits_{k=1}^2d_k\left(q_k^{\beta}+q_k^{N-\beta}\right)\right),\ \beta=1,2,...,N-1,\\
\mathring{C}_{N,0}[a,b]&=&h\left(\frac{1}{2}-\sum\limits_{k=1}^2d_k\frac{q_k-q_k^N}{1-q_k}\right),
\end{eqnarray*}
where $d_k$ are defined from the system
\begin{equation}\label{get_d}
\begin{aligned}
 \sum\limits_{k=1}^2d_k\frac{q_k+q_k^{N+1}}{(q_k-1)^2}&=& \frac{1}{12},\\
\sum\limits_{k=1}^2d_k\frac{q_k-q_k^{N+2}}{(q_k-1)^3}&=& \frac{-1}{24},
\end{aligned}
\end{equation}
for $\omega\in \mathbb{R}$ and $\omega h \not\in \mathbb{Z}$ are expressed as
\begin{eqnarray*}
\mathring{C}_{0,\omega}[a,b]&=&h\left(\frac{\E^{2\pi \i\omega (a+h)}K_{\omega,3}}{\E^{2\pi \i\omega h}-1}-\frac{\E^{2\pi \i\omega a}}{2\pi \i\omega h}+\sum\limits_{k=1}^2\left(a_k\frac{q_k}{q_k-1}+b_k\frac{q_k^N}{1-q_k}\right)\right),   \\
\mathring{C}_{\beta,\omega}[a,b]&=&h\left(\E^{2\pi \i\omega (h\beta+a)}K_{\omega,3}+
\sum\limits_{k=1}^2\left(a_kq_k^\beta+b_kq_k^{N-\beta}\right)\right),
\ \beta=1,2,...,N-1,  \\
\mathring{C}_{N,\omega}[a,b]&=&h\left(\frac{\E^{2\pi \i\omega b}K_{\omega,3}}{1-\E^{2\pi \i\omega h}}+\frac{\E^{2\pi \i\omega b}}{2\pi \i\omega h}+\sum\limits_{k=1}^2\left(a_k\frac{q_k^N}{1-q_k}+b_k\frac{q_k}{q_k-1}\right)\right),
\end{eqnarray*}
where $a_k$ and $b_k$ are determined from the system
\begingroup\makeatletter\def\f@size{8}\check@mathfonts
\begin{equation}\label{get_ab1}
\left.\begin{aligned}
\sum\limits_{k=1}^2a_k\left[\frac{q_k}{(q_k-1)^2}\right]+\sum\limits_{k=1}^2b_k\left[\frac{q_k^{N+1}}{(1-q_k)^2}\right]&=
\frac{\E^{2\pi\i\omega a}}{(2\pi\i\omega h)^2}-\frac{\E^{2\pi\i\omega (h+a)}K_{\omega,3}}{(\E^{2\pi\i\omega h}-1)^2},\\
\sum\limits_{k=1}^2a_k\left[\frac{q_k}{(q_k-1)^3}\right]+\sum\limits_{k=1}^2b_k\left[\frac{q_k^{N+2}}{(1-q_k)^3}\right]&
=\frac{\E^{2\pi\i\omega a}}{(2\pi\i\omega h)^3}-\frac{\E^{2\pi\i\omega a}}{2(2\pi\i\omega h)^2}-\frac{\E^{2\pi\i\omega (h+a)}K_{\omega,3}}{(\E^{2\pi\i\omega h}-1)^3},\\
\sum\limits_{k=1}^2a_k\left[\frac{q_k^{N+1}}{(1-q_k)^2}\right]+\sum\limits_{k=1}^2b_k\left[\frac{q_k}{(q_k-1)^2}\right]&=
\frac{\E^{2\pi\i\omega b}}{(2\pi\i\omega h)^2}-\frac{\E^{2\pi\i\omega (h+b)}K_{\omega,3}}{(1-\E^{2\pi\i\omega h})^2},\\
\sum\limits_{k=1}^2a_k\left[\frac{q_k^2-q_k^{N+2}}{(1-q_k)^3}\right]+\sum\limits_{k=1}^2b_k\left[\frac{q_k^{N+1}-q_k}{(q_k-1)^3}\right]&\\
=(\E^{2\pi\i\omega b}-\E^{2\pi\i\omega a})&\left(\frac{1}{(2\pi\i\omega h)^3}  +\frac{1}{2(2\pi\i\omega h)^2}+\frac{\E^{4\pi\i\omega h}K_{\omega,3}}
{(1-\E^{2\pi\i\omega h})^3}\right),\\
\end{aligned}\right.
\end{equation}
\endgroup
and for $\omega h\in \mathbb{Z}$ with $\omega \neq 0$ have the form
\begin{eqnarray*}
\mathring{C}_{0,\omega}[a,b]&=&h\left(-\frac{\E^{2\pi \i\omega a}}{2\pi \i\omega h}+\sum\limits_{k=1}^2\left(a_k\frac{q_k}{q_k-1}+b_k\frac{q_k^N}{1-q_k}\right)\right),  \\
\mathring{C}_{\beta,\omega}[a,b]&=&h\left(\sum\limits_{k=1}^2\left(a_kq_k^\beta+b_kq_k^{N-\beta}\right)\right),
\ \beta=1,2,...,N-1,  \\
\mathring{C}_{N,\omega}[a,b]&=&h\left(\frac{\E^{2\pi \i\omega b}}{2\pi \i\omega h}+\sum\limits_{k=1}^2\left(a_k\frac{q_k^N}{1-q_k}+b_k\frac{q_k}{q_k-1}\right)\right),
\end{eqnarray*}
where $a_k$ and $b_k$ are defined as follows
\begingroup\makeatletter\def\f@size{10}\check@mathfonts
\begin{equation}\label{get_ab2}
\begin{aligned}
\sum\limits_{k=1}^2a_k\left[\frac{q_k}{(q_k-1)^2}\right]+\sum\limits_{k=1}^2b_k\left[\frac{q_k^{N+1}}{(1-q_k)^2}\right]&=\frac{\E^{2\pi\i\omega a}}{(2\pi\i\omega h)^2},\\
\sum\limits_{k=1}^2a_k\left[\frac{q_k}{(q_k-1)^3}\right]+\sum\limits_{k=1}^2b_k\left[\frac{q_k^{N+2}}{(1-q_k)^3}\right]&=\frac{\E^{2\pi\i\omega a}}{(2\pi\i\omega h)^3}-\frac{\E^{2\pi\i\omega a}}{2(2\pi\i\omega h)^2},\\
\sum\limits_{k=1}^2a_k\left[\frac{q_k^{N+1}}{(1-q_k)^2}\right]+\sum\limits_{k=1}^2b_k\left[\frac{q_k}{(q_k-1)^2}\right]&=\frac{\E^{2\pi\i\omega b}}{(2\pi\i\omega h)^2},\\
\sum\limits_{k=1}^2a_k\left[\frac{q_k^2-q_k^{N+2}}{(1-q_k)^3}\right]+\sum\limits_{k=1}^2b_k\left[\frac{q_k^{N+1}-q_k}{(q_k-1)^3}\right]&\\
=(\E^{2\pi\i\omega b}-\E^{2\pi\i\omega a}&)\left(\frac{1}{(2\pi\i\omega h)^3}+\frac{1}{2(2\pi\i\omega h)^2}\right),\\
\end{aligned}
\end{equation}
\endgroup
$q_k$, $k=1,2$, are roots of the Euler-Frobenius polynomial
\begin{equation}\label{eqE}
E_4(x)=x^4+26x^3+66x^2+26x+1
\end{equation}
with $|q_k|<1$, and
\begin{equation}\label{Kw3}
K_{\omega,3}=\left(\frac{\sin\pi\omega h}{\pi\omega h}\right)^6\frac{5!}{2[\cos4\pi\omega h+26\cos 2\pi\omega h]+66}.
\end{equation}
\end{corollary}

\section{{Application: CT image reconstruction}}

{In this section, we give numerical results of CT image reconstruction by applying the second and third-order optimal quadrature formulas (Corollaries \ref{Cor2} and \ref{Cor3}) to calculate the Fourier} {transform} {and its inversion.
One of the most commonly used CT reconstruction algorithms is FBP {\cite{Buzug08,Feeman15,KakSlaney88,Nat2001}} and as in \cite{KakSlaney88}, the implementation of the FBP consists of four steps: (1) sinogram acquisition, (2) Fourier transform of the sinogram, (3) application of {R}am-{L}ak filter and the Fourier inversion, and (4) back-projection.
To show the effect of the proposed optimal quadrature formula, we implement the FBP in two different ways: one uses \textit{fft} and \textit{ifft} for the Fourier {transform} and its inversion, respectively, and the other uses the proposed optimal quadrature formulas.
{T}hen, we compare the quality of the resulting reconstructed images.
Also, we test the performance of the proposed algorithm by adding the {10\% of} Poisson noise to the given sinogram.
In this section, we deal with the second and third order optimal quadrature formula{s only}.
For the numerical results with the first order optimal quadrature formula, see \cite{HayJeonLee19}.}

\begin{algorithm}
\caption{Reconstruction algorithm with optimal quadrature formula} \label{pseudo_oqf}
\begin{algorithmic}[1]
\State{{{A sinogram {$P(t_m,\theta_k)\,\, \textrm{for} \,\, t_m \in [a, b], \,\, \theta\in [0, \pi]$ is given} as a discrete form.}}}
\State{{{Compute the Fourier transform using the proposed optimal quadrature:}}}
$$S(\omega,\theta)\cong S(\omega,\theta_k)= \sum_{m=0}^{M}\mathring{C}_{m,-\omega}P(t_{m},\theta_{k}), \quad \omega\in \mathbb{R}.$$
\State{{{Compute the inverse Fourier transform using the proposed optimal quadrature:}}}
$$Q(t,\theta)\cong Q(t,\theta_{k}) = \sum_{n=0}^{N}\mathring{C}_{n,t} S(\omega_{n},\theta_{k}) |\omega_{n}|.$$
\State{{{Reconstruct the CT image using back-projection:}}}
$$f(x, y)=\int_{0}^{\pi} Q(t,\theta) d\theta \cong \frac{\pi}{K} \sum_{k=0}^{K-1} Q(t,\theta_{k}).$$
\end{algorithmic}
\end{algorithm}


\begin{figure}
\centering
\includegraphics[width=10cm]{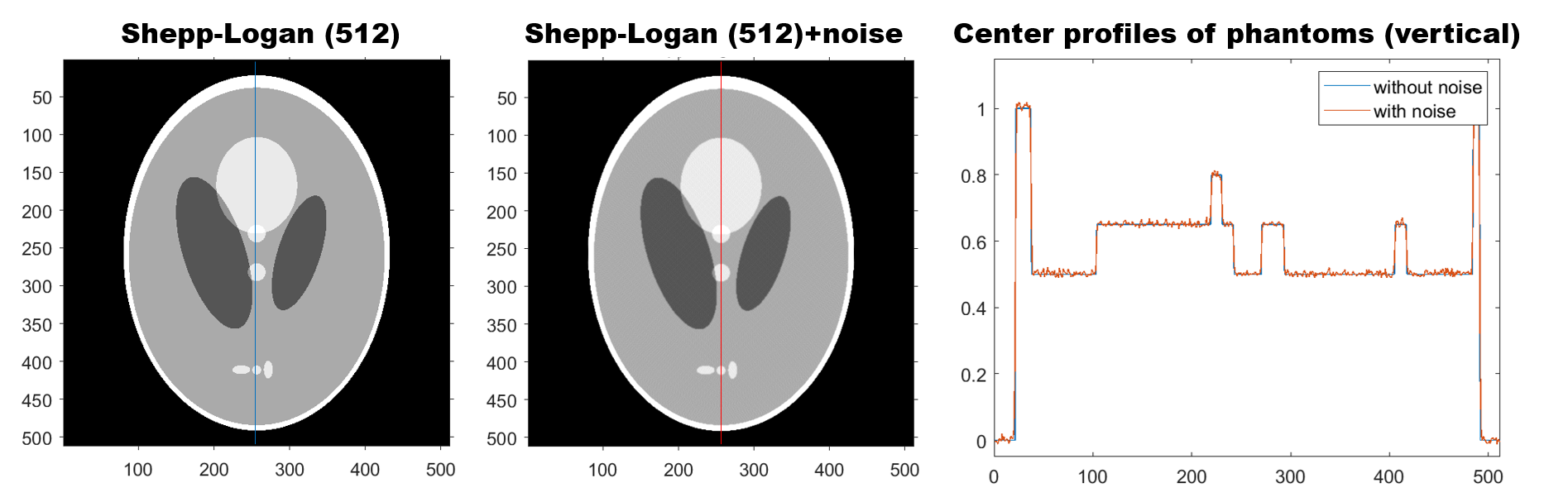}
\caption{{Shepp-Logan phantom images of size $512\times 512$: without noise (left) and with noise (right)} \label{original}}
\end{figure}

Algorithm~\ref{pseudo_oqf} is the pseudo code of the algorithm for CT image reconstruction using optimal quadrature formulas.
For \textit{Step 2} and \textit{Step 3}, the optimal quadrature formulas are used for approximating Fourier integrals.
For the numerical experiment, we use {the Shepp-Logan phantom (SL) {in \cite{KakSlaney88}} and we also consider the {case with noise}} (see Fig.~\ref{original}).
Both are of size {$512 \times 512$} and the sinograms are generated using half rotation sampling with sampling angle~{$0.5^{\circ}$}.
For the implementation, we use $q_{k}=\frac{t_k +\sqrt{t_k^2-4}}{2} (k=1, 2)$, where $t_{k}=-13\pm\sqrt{105}$, the roots of \eqref{eqE} with $|q_k|<1$.
For $q_1$ and $q_2$, $d_k$ are obtained from the $2\times 2$ linear system \eqref{get_d}, and $a_k$ and $b_k$ are obtained from the $4\times 4$ linear systems \eqref{get_ab1} and \eqref{get_ab2}.
{For the numerical experiments, MATLAB 2019b is used.}
For the image quality analysis, we compare maximum error~($E_{\max}$), mean squared error (MSE), and the peak signal-to-noise ratio~(PSNR):~
\begin{eqnarray*}
E_{\max}(I)&=&\max_{i,j}|I(i,j)-I_{ref}(i,j)|, \\
\mbox{MSE}(I)&=&\frac{1}{mn}\sum_{i=1}^{m}\sum_{j=1}^{n}|I(i,j)-I_{ref}(i,j)|^2\,, \\
\mbox{PSNR}(I)&=&10\log_{10}   \left( \frac{I_{\max}^2}{\mbox{MSE}(I)} \right)\,,
\end{eqnarray*}
where $I_{\max}$ is the maximum pixel value of the image $I$.
Images of original simulated phantoms are denoted by $I_{ref}$ (Fig.~\ref{original}).


\begin{figure}
\centering
\includegraphics[width=10cm]{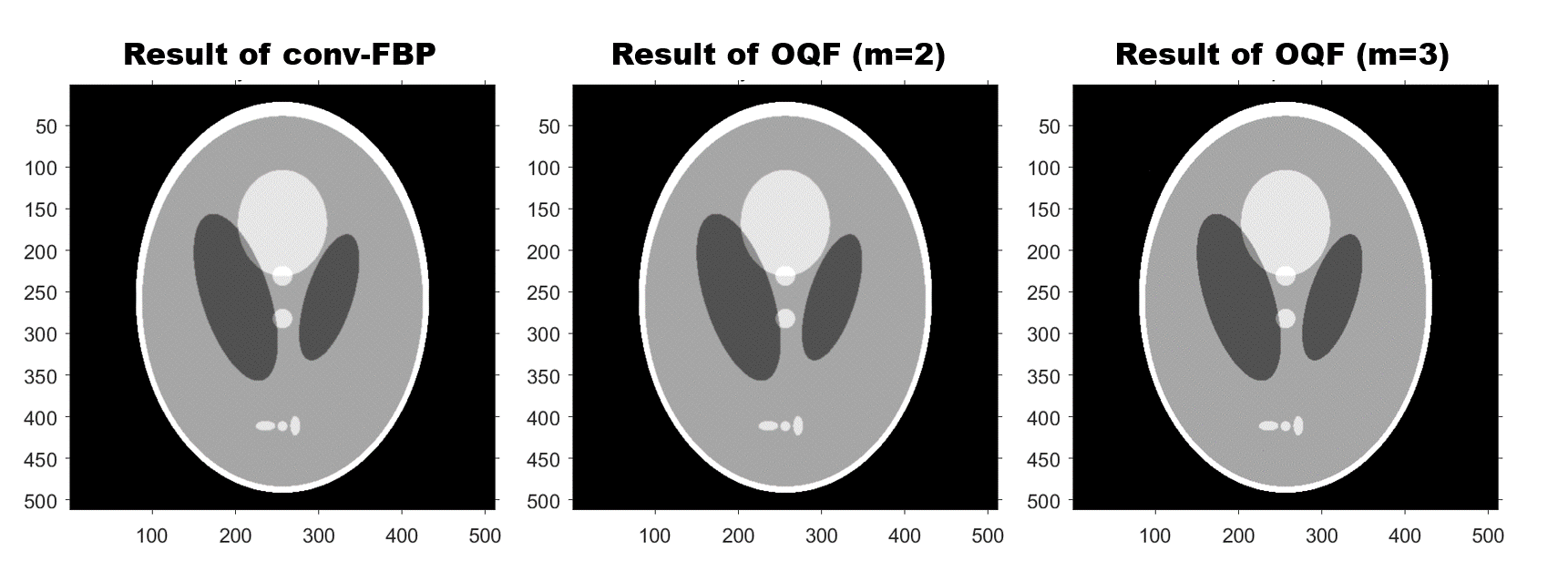}
\caption{Reconstruction results of the conventional FBP using \textit{fft} and \textit{ifft} (left), using the second order optimal quadrature formula (middle), and using the third order optimal quadrature formula (right). \label{recon_nonoise}}
\end{figure}

\begin{figure}
\centering
\includegraphics[width=10cm]{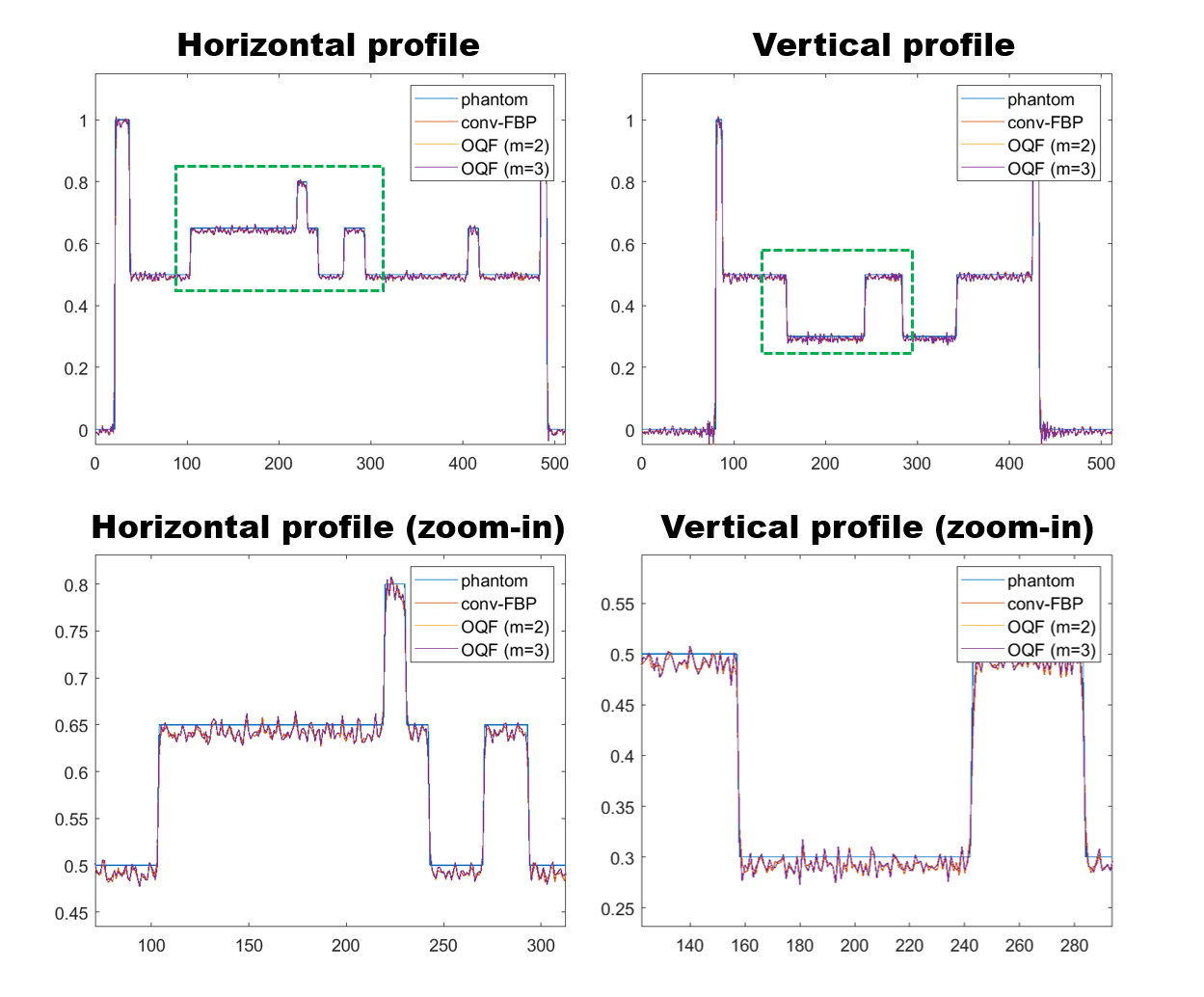}
\caption{Profiles of the reconstructed image of the Shepp-Logan: horizontal (left) and vertical (right). \label{profile}}
\end{figure}

\begin{figure}
\centering
\includegraphics[width=10cm]{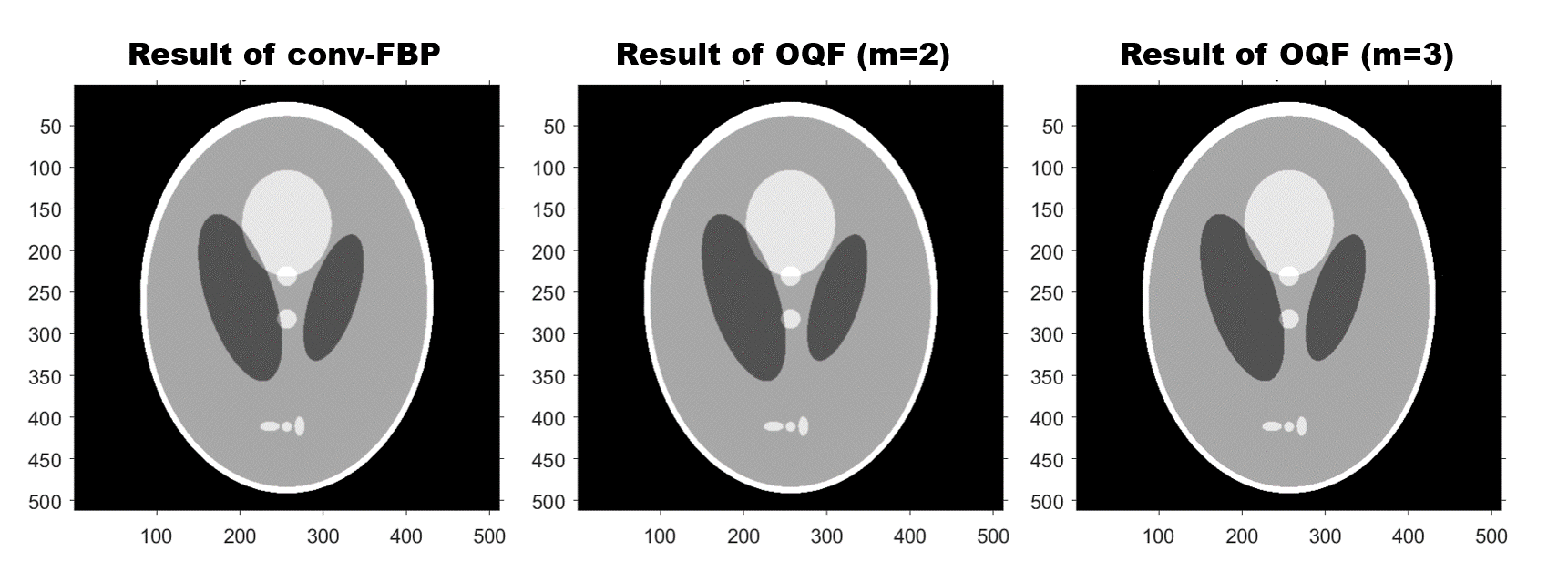}
\caption{\textbf{Noise test}: reconstruction results of the conventional FBP using \textit{fft} and \textit{ifft} (left), using the second order optimal quadrature formula (middle), and using the third order optimal quadrature formula (right). \label{recon_noise}}
\end{figure}

\begin{figure}
\centering
\includegraphics[width=10cm]{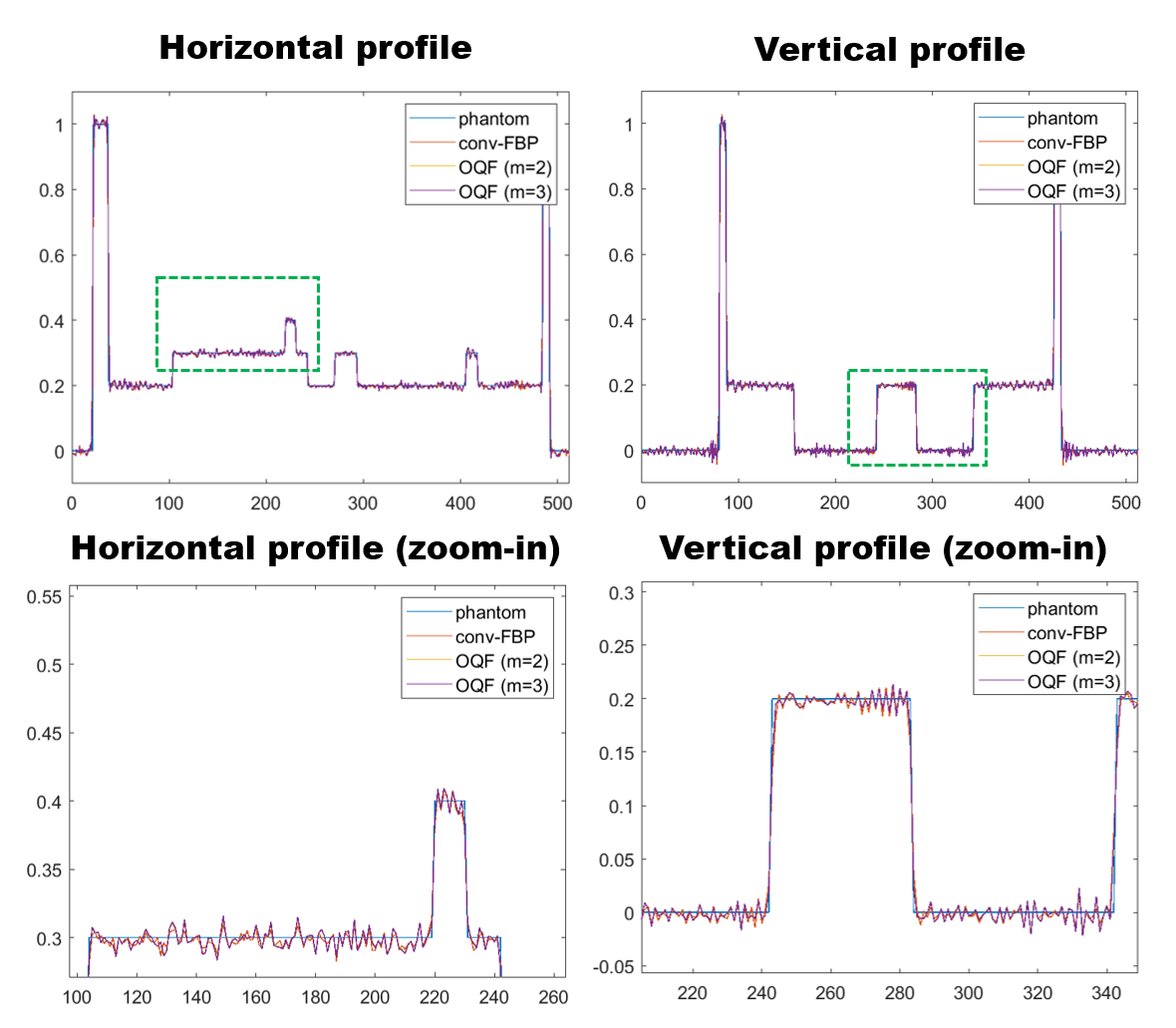}
\caption{\textbf{Noise test}: Profiles of the reconstructed image of the Shepp-Logan: horizontal (left) and vertical (right). \label{profile_noise}}
\end{figure}

{
\begin{table}[]
\centering
\begin{tabular}{|c||r|r|r||r|r|r|}
\hline
\multicolumn{1}{|l||}{\multirow{2}{*}{}}       & \multicolumn{3}{c||}{Shepp-Logan without noise}                            & \multicolumn{3}{c|}{Shepp-Logan with noise}                            \\
 \cline{2-7}
\multicolumn{1}{|c||}{}                  & \multicolumn{1}{c|}{\begin{tabular}[c]{@{}c@{}} conv-FBP\\ \textit{fft} and \textit{ifft}\end{tabular}} & \multicolumn{1}{c|}{\begin{tabular}[c]{@{}c@{}}OQF\\ $m=2$\end{tabular}}& \multicolumn{1}{c||}{\begin{tabular}[c]{@{}c@{}}OQF\\ $m=3$\end{tabular}} &   \multicolumn{1}{c|}{\begin{tabular}[c]{@{}c@{}} conv-FBP\\ \textit{fft} and \textit{ifft}\end{tabular}} & \multicolumn{1}{c|}{\begin{tabular}[c]{@{}c@{}}OQF\\ $m=2$\end{tabular}}& \multicolumn{1}{c|}{\begin{tabular}[c]{@{}c@{}}OQF\\ $m=3$\end{tabular}}  \\ \hline \hline
$E_{\max}$   &   0.3458    &  0.3526    &   0.3307    &  0.3722      & 0.3634       &  0.3472           \\ \hline
MSE   & 7.9648e-04  &  7.2111e-04   &    6.5084e-04     &7.9088e-04    &7.4509e-04  &   6.4990e-04     \\ \hline
PSNR   &    30.9883   &    31.4200  &  31.8652    &  31.0189     & 31.2779&  31.8715    \\ \hline
\end{tabular}
\caption{{Quantitative analysis for the reconstructed CT image from FBP using conventional \textit{fft-ifft} and the second and third order optimal quadrature formula{s} in Corollaries \ref{Cor2} and \ref{Cor3}. }}
\label{imageAnalysis_23}
\end{table}}

Figs.~\ref{recon_nonoise} and \ref{profile} show the results of CT image reconstruction using optimal quadrature formulas of the second and third orders and profiles of the results, respectively.
As shown in Fig.~\ref{profile}, the profile lines of {conv-FBP} and optimal quadrature formula almost coincide.
{Figs.~\ref{recon_noise} and \ref{profile_noise} show the results of noise test and we observe that the resulting reconstructed images are very similar, visually.}
Table~\ref{imageAnalysis_23} shows $E_{\max}$, MSE, and PSNR for the reconstruction results using the second and third order optimal quadrature formulas.
{The second order optimal quadrature formula performs similar to or slightly better than the conventional FBP.
The proposed third order optimal quadrature formula produce{s} more improved quality than the conventional FBP particularly in terms of MSE and PSNR.
The third order optimal quadrature formula also show{s} the best performance when noise was added.}

\textit{\bf{Remark 3.}} The computational complexity of optimal quadrature formulas for Fourier integration is $\mathcal{O}(n^2)$ for the $n\times n$ matrix multiplication and addition and it is larger than $\mathcal{O}(n\log n)$ for the fast Fourier transformation.
However, once the coefficients of optimal quadrature formulas for Fourier integral are computed, they can be used several times for the same size of input signals. Moreover, they can be well optimized using parallelization.

\textit{\bf{Remark 4.}} It should be noted that since the coefficients of the optimal quadrature formulas of the form (\ref{eq2.28}) are continuous functions of the parameter $\omega$, it is not needed any interpolation in
Algorithm 1.

\section*{Acknowledgements}
The work has been done while A.R. Hayotov was visiting Department of Mathematical
Sciences at KAIST, Daejeon, Republic of Korea.
A.R. Hayotov's work was supported by the `Korea Foundation for Advanced Studies'/`Chey Institute for Advanced Studies' International Scholar Exchange Fellowship for academic year of 2018-2019. C.-O. Lee's work was supported by {the National Research Foundation of Korea (NRF)} grant funded by {the Korea government (MSIT)} (No.~NRF-2017R1A2B4011627).


\end{document}